\documentclass[final,1p,times]{elsarticle}

\usepackage{amssymb}
 \usepackage{amsthm}
\usepackage{amscd}
\usepackage{amsmath}
\usepackage{amsfonts}
\usepackage{amssymb}
\usepackage{graphicx}
\newtheorem{theorem}{Theorem}

\usepackage{mathrsfs}
\usepackage{titletoc}

\newcommand{\ra}{\rightarrow}
\newcommand{\p}{\partial}
\newcommand{\f}{\frac}

\newcommand{\be}{\begin{equation}}
\renewcommand{\ra}{\rightarrow}
\newcommand{\ee}{\end{equation}}
\newcommand{\bea}{\begin{eqnarray}}
\newcommand{\eea}{\end{eqnarray}}
\newcommand{\bna}{\begin{eqnarray*}}
\newcommand{\ena}{\end{eqnarray*}}

\renewcommand{\le}{\left}
\newcommand{\ri}{\right}

\journal{***}

\begin{document}

\begin{frontmatter}
\title{Yamabe type equations on graphs}

\author{Alexander Grigor'yan}
\ead{grigor@math.uni-bielefeld.de}
\address{Department of Mathematics,
University of Bielefeld, Bielefeld 33501, Germany}

\author{Yong Lin}
 \ead{linyong01@ruc.edu.cn}

\author{Yunyan Yang}
 \ead{yunyanyang@ruc.edu.cn}

\address{ Department of Mathematics,
Renmin University of China, Beijing 100872, P. R. China}

\begin{abstract}
Let $G=(V,E)$ be a locally finite graph, $\Omega\subset V$ be a bounded domain, $\Delta$ be the usual
graph Laplacian, and $\lambda_1(\Omega)$ be the first eigenvalue of $-\Delta$ with respect to Dirichlet boundary
condition. Using the mountain pass theorem due to Ambrosetti-Rabinowitz, we prove that if $\alpha<\lambda_1(\Omega)$, then for any $p>2$,
there exists a positive solution to
$$\le\{\begin{array}{lll}
-\Delta u-\alpha u=|u|^{p-2}u\,\,\,{\rm in}\,\,\, \Omega^\circ,\\[1.5ex]
 u=0\,\,\,{\rm on}\,\,\, \p\Omega,
\end{array}\ri.
$$
where $\Omega^\circ$ and $\p\Omega$ denote the interior and the boundary of $\Omega$
 respectively.
 Also we consider similar problems involving the $p$-Laplacian and poly-Laplacian by the same method.
 Such problems can be viewed as discrete versions of  the Yamabe type equations on Euclidean space or compact Riemannian manifolds.
\end{abstract}


\end{frontmatter}

\titlecontents{section}[0mm]
                       {\vspace{.2\baselineskip}}
                       {\thecontentslabel~\hspace{.5em}}
                        {}
                        {\dotfill\contentspage[{\makebox[0pt][r]{\thecontentspage}}]}
\titlecontents{subsection}[3mm]
                       {\vspace{.2\baselineskip}}
                       {\thecontentslabel~\hspace{.5em}}
                        {}
                       {\dotfill\contentspage[{\makebox[0pt][r]{\thecontentspage}}]}

\setcounter{tocdepth}{2}


\section{Introduction}

Let $\Omega$ be a domain of $\mathbb{R}^n$ and $W_0^{1,q}(\Omega)$ be the closure of $C_0^\infty(\Omega)$ with respect to the norm
$$\|u\|_{W_0^{1,q}(\Omega)}=\le(\int_\Omega\le(|\nabla u|^q+|u|^q\ri)dx\ri)^{1/q},$$
where $q\geq 1$. Then the Sobolev embedding theorem reads
$$W_0^{1,q}(\Omega)\hookrightarrow\le\{\begin{array}{lll}
L^p(\Omega)\,\,\,{\rm for}\,\,\, q\leq p\leq q^*=\f{nq}{n-q},&{\rm when}\,\,\, n>q\\[1.5ex]
L^p(\Omega)\,\,\,{\rm for}\,\,\, q\leq p<+\infty, &{\rm when}\,\,\, q=n\\[1.5ex]
C^{1-n/q}(\overline{\Omega}), &{\rm when}\,\,\, q>n.
\end{array}\ri.$$
 The model problem
\be\label{P}\le\{\begin{array}{lll}
-\Delta u+\lambda u=|u|^{p-2}u,\\[1.5ex]
u\in W_0^{1,2}(\Omega)
\end{array}\ri.\ee
or its variants has been extensively studied since 1960's. Let $J:W_0^{1,2}(\Omega)\ra \mathbb{R}$
be a functional defined by
$$J(u)=\f{1}{2}\int_\Omega(|\nabla u|^2+\lambda u^2)dx-\f{1}{p}\int_\Omega|u|^pdx.$$
Clearly the critical points of $J$ are weak solutions to the problem (\ref{P}). In the case $2<p<+\infty$ when $n=1,\,2$, or
$2<p\leq 2^*=2n/(n-2)$ when $n\geq 3$, one can check that $\sup_{W_0^{1,2}(\Omega)}J=+\infty$ and $\inf_{W_0^{1,2}(\Omega)}J
=-\infty$. In \cite{Nehari-1}, Nehari  obtained a nontrivial solution of (\ref{P}) when $\lambda\geq 0$ and $\Omega=(a,b)$, by minimizing
$J$ in the manifold $$\mathcal{N}=\{u\in W_0^{1,2}(\Omega):\langle J^\prime(u),u\rangle=0,\,u\not= 0\}.$$ In \cite{Nehari-2},
he proved the existence of infinitely many solutions  and in \cite{Nehari-3}, he solved the case
where $\Omega=\mathbb{R}^3$, $\lambda>0$ and $2<p<6$ after reduction to an ordinary differential equation. When
$\Omega$ is unbounded or when $p=2^*$, there is a lack of compactness in Sobolev spaces
because of invariance by translation or
  by dilation. Some nonexistence results follow from the Pohozaev identity. General existence theorems were first
obtained by Strauss \cite{Strauss} when $\Omega=\mathbb{R}^n$ and by Brezis-Nirenberg \cite{Brezis-Nirenberg} when $p=2^*$. When $p=2^*$, Bahri and Coron \cite{Bahri-Coron} proved that if there exist a positive integer $d$ such that the $d$-dimentional homology group of domain $\Omega$ is nontrivial, then (\ref{P}) has a solution. Pohozaev
\cite{Pohozaev1} proved that if $\Omega$ is starshaped, then (\ref{P}) has no solution.
The Brezis-Lieb lemma and Lions' concentration compactness principle are important tools in solving those problems. For other existence results
for variants of (\ref{P}), we refer the reader to \cite{Willem}.

Analogous to (\ref{P}), one can consider the problem
\be\label{PT}\le\{\begin{array}{lll}
-\Delta_n u=f(x,u) \quad {\rm in}\quad \Omega\\[1.5ex]
u\in W_0^{1,n}(\Omega),
\end{array}\ri.\ee
where $\Delta_n$ is the $n$-Laplace operator and $f(x,s)$ has exponential growth as $s\ra+\infty$. Instead of the Sobolev embedding
theorem, the key tool in solving the problem (\ref{PT}) is the Trudinger-Moser embedding contributed by Yudovich \cite{Yudovich},
Pohozaev \cite{Pohozaev}, Peetre \cite{Peetre}, Trudinger \cite{Trudinger} and Moser \cite{Moser}. In \cite{Adimurthi}, Adimurthi
proved an existence
of positive solution to (\ref{PT}) by using a method of Nehari manifold.
In \cite{FMR}, de Figueiredo,  Miyagaki and Ruf  considered (\ref{PT}) in the case that $\Omega$ is a bounded domain in $\mathbb{R}^2$,
by using the critical point theory. In \cite{doo1}, by using the mountain-pass theorem
without the Palais-Smale condition, do \'O  improve the results of \cite{Adimurthi,FMR}.
In \cite{doo2}, using the same method, he extended these results to the case that $\Omega$ is the whole Euclidean space $\mathbb{R}^n$.
For related works, we refer the reader to \cite{doo3,Adi-Yang,Yang1,Yang2} and the references there in.

On Riemannian manifolds, an analog of the model problem (\ref{P}) arises from the Yamabe problem: Let $(M,g)$ be a
compact $n$ $(\geq 3)$ dimensional Riemannian manifold without boundary. Does there exist a good metric $\tilde{g}$ in the conformal class of $g$ such that
the scalar curvature $R_{\tilde{g}}$ is a constant? This problem was studied by Yamabe \cite{Yamabe}, Trudinger \cite{Tru},
Aubin \cite{Aubin}, and completely solved by Schoen \cite{Schoen}. Though there is no background of geometry or physics,
there are still some works concerning the problem (\ref{PT})
on Riemannian manifolds, see for examples \cite{Yang-Zhao,do-Yang,Zhao,Yang-jfa}.\\

Our goal is to consider problems (\ref{P}) and (\ref{PT}) when
an Euclidean domain $\Omega$ is replaced by a graph. Such problems can be viewed as discrete versions of (\ref{P}) and (\ref{PT}).
In this paper, we concern bounded domain on locally finite graphs or finite graphs.
The key point is an observation of pre-compactness of the Sobolev
space in our setting.
Using the mountain pass theorem due to Ambrosseti-Rabinowich \cite{Ambrosetti-Rabinowitz}, we
 prove the existence of nontrivial solutions to Yamabe type equations on graphs. Our results are quite different from
 that of  \cite{Bahri-Coron} and \cite{Pohozaev1} when $p=2^*$, namely, the existence results are not related to the homology group of the graphs.\\

 This paper is organized as follows: In Section 2, we give some notations on graph and state main results.
 In Section 3, we establish Sobolev embedding such that the mountain pass theorem can be applied to our problems.
 Local existence results (Theorems \ref{Theorem 1}-\ref{Theorem 4}) are proved in Section 4, and
 global existence results (Theorems \ref{Theorem 5}-\ref{Theorem 8})
 are proved in Section 5.

\section{Settings and main results}

Let $G=(V,E)$ be a finite or locally finite graph, where $V$ denotes the vertex set and $E$ denotes the edge set.
For any edge $xy\in E$, we assume that its weight $w_{xy}>0$ and that $w_{xy}=w_{yx}$. The degree of $x\in V$ is defined
as ${\rm deg}(x)=\sum_{y\sim x}w_{xy}$, where we write $y\sim x$ if
$xy\in E$. Let $\mu:V\ra \mathbb{R}^+$ be a finite measure. For any function $u:V\ra \mathbb{R}$, the $\mu$-Laplacian (or Laplacian for short)
of $u$ is defined as
\be\label{lap}\Delta u(x)=\f{1}{\mu(x)}\sum_{y\sim x}w_{xy}(u(y)-u(x)).\ee
The associated gradient form reads
\be\label{grad-form}\Gamma(u,v)(x)=\f{1}{2\mu(x)}\sum_{y\sim x}w_{xy}(u(y)-u(x))(v(y)-v(x)).\ee
Write $\Gamma(u)=\Gamma(u,u)$. We denote the length of its gradient by
\be\label{grd}|\nabla u|(x)=\sqrt{\Gamma(u)(x)}=\le(\f{1}{2\mu(x)}\sum_{y\sim x}w_{xy}(u(y)-u(x))^2\ri)^{1/2}.\ee
Similar to the Euclidean case, we define the length of $m$-order gradient of $u$ by
 $$\label{high-deriv}|\nabla^mu|=\le\{\begin{array}{lll}
 |\nabla\Delta^{\f{m-1}{2}}u|,\,\,{\rm when}\,\,\, m
 \,\,{\rm is\,\,odd}\\[1.5ex]
 |\Delta^{\f{m}{2}}u|,\,\,{\rm when}\,\,\, m
 \,\,{\rm is\,\,even},
 \end{array}\ri.
 $$
 where $|\nabla\Delta^{\f{m-1}{2}}u|$ is defined as in (\ref{grd}) with $u$ is replaced by $\Delta^{\f{m-1}{2}}u$, and
 $|\Delta^{\f{m}{2}}u|$ denotes the usual absolute of the function $\Delta^{\f{m}{2}}u$.
 Let $\Omega$ be a domain (or a connected set) in $V$.  
 To compare with the Euclidean setting, we denote, for any function $u:V\ra\mathbb{R}$,
 \be\label{integration}\int_\Omega ud\mu=\sum_{x\in \Omega}\mu(x)u(x).\ee
 The first eigenvalue of the Laplacian with respect to Dirichlet boundary condition reads
 \be\label{eigen}\lambda_1(\Omega)=\inf_{u\not\equiv 0,\,u|_{\p\Omega}=0}\f{\int_\Omega|\nabla u|^2d\mu}{\int_\Omega u^2d\mu},\ee
 where $\p\Omega$ is the boundary of $\Omega$, namely
 $\p\Omega=\{x\in\Omega: \exists y\not\in \Omega\,\,{\rm such\,\,that}\,\, xy\in E\}$.
 Moreover, we denote the interior of $\Omega$ by
 $\Omega^\circ=\Omega\setminus\p\Omega$. Our first result is the following:

\begin{theorem}\label{Theorem 1}
 Let $G=(V,E)$ be a locally finite graph, $\Omega\subset V$ be a bounded domain with $\Omega^\circ\not=\varnothing$, and
 $\lambda_1(\Omega)$ be defined as in (\ref{eigen}). Then for any $p>2$ and any $\alpha<\lambda_1(\Omega)$, there exists a
 solution to the equation
 $$\label{1}\left\{
 \begin{array}{lll}
 -\Delta u-\alpha u=|u|^{p-2}u\quad{\rm in}\quad \Omega^\circ\\[1.5ex]
 u>0\,\,\,{\rm in}\,\,\, \Omega^\circ,\,\,\, u=0\,\,\,{\rm on}\,\,\, \p\Omega.
 \end{array}
 \right.
 $$
 \end{theorem}

The $p$-Laplacian of $u:V\ra\mathbb{R}$, namely $\Delta_pu$, is defined in the distributional sense by
$$\int_V(\Delta_pu)\phi d\mu=-\int_V |\nabla u|^{p-2}\Gamma(u,\phi)d\mu,\quad\forall
\phi\in \mathcal{C}_{\rm c}(V),$$
 where $\Gamma(u,\phi)$ is defined as in (\ref{grad-form}), $\mathcal{C}_{\rm c}(V)$ denotes
the set of all functions with compact support, and the integration is defined by (\ref{integration}). Point-wisely, $\Delta_pu$ can be written as
$$\Delta_pu(x)=\f{1}{2\mu(x)}\sum_{y\sim x}\le(|\nabla u|^{p-2}(y)+|\nabla u|^{p-2}(x)\ri)w_{xy}(u(y)-u(x)).$$
When $p=2$, $\Delta_p$ is the standard graph Laplacian $\Delta$ defined by (\ref{lap}).
   The first eigenvalue of the $p$-Laplacian with respect to Dirichlet boundary
condition reads
\be\label{p-eig}\lambda_p(\Omega)=\inf_{u\not\equiv 0,\,u|_{\p\Omega}=0}\f{\int_\Omega|\nabla u|^pd\mu}
{\int_\Omega |u|^pd\mu}.\ee
Our second result can be stated as follows:
\begin{theorem}\label{Theorem 2}
 Let $G=(V,E)$ be a locally finite graph, $\Omega\subset V$ be a bounded domain with $\Omega^\circ\not=\varnothing$.
 Let
 $\lambda_p(\Omega)$ be defined as in (\ref{p-eig}) for some $p>1$. Suppose that $f:\Omega\times\mathbb{R}\ra \mathbb{R}$ satisfies the following
 hypothesis:\\
 $(H_1)$ For any $x\in \Omega$, $f(x,t)$ is continuous in $t\in\mathbb{R}$;\\
 $(H_2)$ For all $(x,t)\in\Omega\times[0,+\infty)$, $f(x,t)\geq 0$ , and $f(x,0)=0$ for all $x\in \Omega$;\\
 $(H_3)$ There exists some $q>p$ and $s_0>0$ such that if $s\geq s_0$, then there holds
 $$F(x,s)=\int_0^sf(x,t)dt\leq \f{1}{q}sf(x,s),\quad \forall x\in \Omega;$$
 $(H_4)$ For any $x\in \Omega$, there holds
 $$\limsup_{t\ra 0+}\f{f(x,t)}{t^{p-1}}<\lambda_p(\Omega).$$
 Then there exists a nontrivial solution to the equation
 $$\label{2}\left\{
 \begin{array}{lll}
 -\Delta_p u=f(x,u)\quad{\rm in}\quad \Omega^\circ\\[1.5ex]
  u\geq 0\,\,\,{\rm in}\,\,\, \Omega^\circ,\,\,\, u=0\,\,\,{\rm on}\,\,\, \p\Omega.
 \end{array}
 \right.
 $$
 \end{theorem}

 In Theorem \ref{Theorem 2}, if $p=2$, then
  $|s|^{q-2}s$ $(q> 2)$ satisfies $(H_1)-(H_4)$. Moreover, the nonlinearities in Theorem \ref{Theorem 2}
 include the case of exponential growth as in the problem  (\ref{PT}). For further extension, we
 define an analog of $\lambda_p(\Omega)$ by
 \be\label{l-mp}\lambda_{mp}(\Omega)=\inf_{u\in\mathcal{H}}\f{\int_\Omega|\nabla^mu|^pd\mu}{\int_\Omega|u|^pd\mu},\ee
 where $m$ is any positive integer and $\mathcal{H}$ denotes the set of all functions $u\not\equiv 0$ with $u=|\nabla u|=\cdots=|\nabla^{m-1}u|=0$
 on $\p\Omega$. Then we have the following:

 \begin{theorem}\label{Theorem 4}
 Let $G=(V,E)$ be a locally finite graph and $\Omega\subset V$ be a bounded domain with $\Omega^\circ\not=\varnothing$.
 Let  $m\geq 2$ be an integer, $p>1$, and $\lambda_{mp}(\Omega)$ be defined by (\ref{l-mp}).  Suppose that $f:\Omega\times\mathbb{R}\ra \mathbb{R}$ satisfies the following
 assumptions:\\
 $(A_1)$ $f(x,0)=0$, $f(x,t)$ is continuous with respect to $t\in\mathbb{R}$;\\
 $(A_2)$  $\limsup_{t\ra 0}\f{|f(x,t)|}{|t|^{p-1}}<\lambda_{mp}(\Omega)$;\\
 $(A_3)$ there exists some $q>p$ and $M>0$ such that if $|s|\geq M$, then
 $$0<qF(x,s)\leq sf(x,s),\quad\forall x\in \Omega.$$
 Then there exists a nontrivial
 solution to the equation
 $$\label{e-3}\left\{
 \begin{array}{lll}
 \mathcal{L}_{m,p}\, u=f(x,u)\quad{\rm in}\quad \Omega^\circ\\[1.5ex]
  |\nabla^j u|=0\,\,{\rm on}\,\, \p\Omega,\,0\leq j\leq m-1,
 \end{array}
 \right.
 $$
 where $\mathcal{L}_{m,p}\, u$ is defined as follows: for any $\phi$ with $\phi=|\nabla\phi|=\cdots=
 |\nabla^{m-1}\phi|=0$ on $\p\Omega$, there holds
 $$\int_\Omega(\mathcal{L}_{m,p}\, u)\phi d\mu=\le\{\begin{array}{lll}
 \int_\Omega|\nabla^mu|^{p-2}\Gamma(\Delta^{\f{m-1}{2}}u,\Delta^{\f{m-1}{2}}\phi)d\mu,\,\,\,{\rm when} \,\, m\,\,{\rm is\,\,odd},\\
 [1.5ex] \int_\Omega |\nabla^mu|^{p-2}\Delta^{\f{m}{2}}u\Delta^{\f{m}{2}}\phi d\mu,\,\,\,{\rm when} \,\, m\,\,{\rm is\,\,even}.
 \end{array}\ri.$$
 In particular, if $p=2$, then $\mathcal{L}_{m,p}\, u=(-\Delta)^mu$, the poly-Laplacian of $u$.
 \end{theorem}

 If $G=(V,E)$ is a finite graph, we also have existence results similar to the above theorems. Analogous to Theorem \ref{Theorem 1},
 we state the following:

 \begin{theorem}\label{Theorem 5}
 Let $G=(V,E)$ be a finite graph. Suppose that $p>2$ and $h(x)>0$ for all $x\in V$. Then there exists a
 solution to the equation
 $$\label{e-4}\left\{
 \begin{array}{lll}
 -\Delta u+h u=|u|^{p-2}u\quad{\rm in}\quad V\\[1.5ex]
 u>0\,\,\,{\rm in}\,\,\, V.
 \end{array}
 \right.
 $$
 \end{theorem}

 Similar to Theorem \ref{Theorem 2}, we have

 \begin{theorem}\label{Theorem 6}
 Let $G=(V,E)$ be a finite graph.  Suppose that  $h(x)>0$ for all $x\in V$.
 Suppose that $f:V\times\mathbb{R}\ra \mathbb{R}$ satisfies the following
 hypothesis:\\
 $(H_V^1)$ For any $x\in V$, $f(x,t)$ is continuous in $t\in\mathbb{R}$;\\
 $(H_V^2)$ For all $(x,t)\in V\times[0,+\infty)$, $f(x,t)\geq 0$ , and $f(x,0)=0$ for all $x\in V$;\\
 $(H_V^3)$ There exists some $q>p>1$ and $s_0>0$ such that if $s\geq s_0$, then there holds
 $$F(x,s)=\int_0^sf(x,t)dt\leq \f{1}{q}sf(x,s),\quad \forall x\in V;$$
 $(H_V^4)$ For any $x\in V$, there holds
 $$\limsup_{t\ra 0+}\f{f(x,t)}{t^{p-1}}<\lambda_p(V)=\inf_{u\not\equiv 0}\f{\int_V(|\nabla u|^p+h|u|^p)d\mu}
 {\int_V|u|^pd\mu}.$$
 Then there exists a nontrivial solution to the equation
 $$\label{e-5}\left\{
 \begin{array}{lll}
 -\Delta_p u+h|u|^{p-2}u=f(x,u)\quad{\rm in}\quad V\\[1.5ex]
  u\geq 0\,\,\,{\rm in}\,\,\, V,
 \end{array}
 \right.
 $$
 where $\Delta_pu$ denotes the $p$-Laplacian of $u$.
 \end{theorem}

 Finally we have an analog of Theorem \ref{Theorem 4}, namely

 \begin{theorem}\label{Theorem 8}
 Let $G=(V,E)$ be a finite graph. Let  $m\geq 2$ be an integer and $p>1$. Suppose that
 $h(x)>0$ for all $x\in V$.  Assume $f(x,u)$ satisfies the following assumptions:\\
 $(A_V^1)$ For any $x\in V$, $f(x,0)=0$, $f(x,t)$ is continuous with respect to $t\in\mathbb{R}$;\\
 $(A_V^2)$  $\limsup_{t\ra 0}\f{|f(x,t)|}{|t|^{p-1}}<\lambda_{mp}(V)=\inf_{u\not\equiv 0}
 \f{\int_V(|\nabla^m u|^p+h|u|^p)d\mu}{\int_V|u|^pd\mu}$;\\
 $(A_V^3)$ there exists some $q>p$ and $M>0$ such that if $|s|\geq M$, then
 $$0<qF(x,s)\leq sf(x,s),\quad\forall x\in V.$$
Then there exists a nontrivial solution to
$$\label{e-6}
\mathcal{L}_{m,p}\, u+h|u|^{p-2}u=f(x,u)\quad{\rm in}\quad V,
$$
 where $\mathcal{L}_{m,p}\, u$ is defined in the distributional sense: for any function $\phi$, there holds
 $$\int_V(\mathcal{L}_{m,p}\, u)\phi d\mu=\le\{\begin{array}{lll}
 \int_V|\nabla^mu|^{p-2}\Gamma(\Delta^{\f{m-1}{2}}u,\Delta^{\f{m-1}{2}}\phi)d\mu,\,\,\,{\rm when} \,\, m\,\,{\rm is\,\,odd},\\
 [1.5ex] \int_V |\nabla^mu|^{p-2}\Delta^{\f{m}{2}}u\Delta^{\f{m}{2}}\phi d\mu,\,\,\,{\rm when} \,\, m\,\,{\rm is\,\,even}.
 \end{array}\ri.$$
  \end{theorem}

\section{Preliminary analysis}

  Let $G=(V,E)$ be a locally finite graph, $\Omega\subset V$ be an domain, $\p\Omega$ be its boundary and $\Omega^\circ$ be its
  interior. For any $p>1$,
 $W^{m,p}(\Omega)$ is defined as a space of all functions $u:V\ra\mathbb{R}$ satisfying
 \be\label{norm}\|u\|_{W^{m,p}(\Omega)}=\le(\sum_{k=0}^m\int_\Omega|\nabla^ku|^pd\mu\ri)^{1/p}<\infty.\ee
 Denote $\mathcal{C}_0^m(\Omega)$  be a set of all functions $u:\Omega\ra\mathbb{R}$ with
 $u=|\nabla u|=\cdots=|\nabla^{m-1}u|=0$ on $\p\Omega$.
 We denote $W_0^{m,p}(\Omega)$ be the completion of $\mathcal{C}_0^m(\Omega)$ under the norm (\ref{norm}).
 If we further assume that $\Omega$ is a bounded, then $\Omega$ is a finite set.
 Observing that the dimension of $W_0^{m,p}(\Omega)$ is finite when $\Omega$ is bounded,
 we have the following Sobolev embedding:

 \begin{theorem} \label{compactlemma} Let $G=(V,E)$ be a locally finite graph,
 $\Omega$ be a bounded domain of $V$ such that $\Omega^\circ\not=\varnothing$.
 Let $m$ be any positive integer and $p>1$.
 Then $W_0^{m,p}(\Omega)$ is embedded in $L^q(\Omega)$ for all $1\leq q\leq+\infty$. In particular,
   there exists a constant $C$ depending only on $m$, $p$ and $\Omega$ such that
   \be\label{Poincare}\le(\int_\Omega |u|^qd\mu\ri)^{1/q}\leq C\le(\int_\Omega|\nabla^mu|^pd\mu\ri)^{1/p}\ee
   for all $1\leq q\leq +\infty$ and for all $u\in W_0^{m,p}(\Omega)$.
  Moreover, $W_0^{m,p}(\Omega)$ is pre-compact, namely, if $u_k$ is bounded in $W_0^{m,p}(\Omega)$, then up to a subsequence,
  there exists some $u\in W_0^{m,p}(\Omega)$ such that $u_k\ra u$ in $W_0^{m,p}(\Omega)$.
 \end{theorem}

 {\it Proof.} Since $\Omega$ is a finite set, $W_0^{m,p}(\Omega)$ is a finite dimensional space. Hence $W_0^{m,p}(\Omega)$
 is pre-compact. We are left to show (\ref{Poincare}). It is not difficult to see that
 \be\label{w0-norm}\|u\|_{W_0^{m,p}(\Omega)}=\le(\int_\Omega|\nabla^mu|^pd\mu\ri)^{1/p}\ee
 is a  norm equivalent to (\ref{norm}) on $W_0^{m,p}(\Omega)$. Hence for any $u\in W_0^{m,p}(\Omega)$, there exists some
 constant $C$ depending only on $m$, $p$ and $\Omega$ such that
 \be\label{qp}\le(\int_\Omega |u|^pd\mu\ri)^{1/p}=\le(\sum_{x\in\Omega}\mu(x)|u(x)|^p\ri)^{1/p}\leq
 C\le(\int_\Omega|\nabla^mu|^pd\mu\ri)^{1/p},\ee
 since $W_0^{m,p}(\Omega)$ is a finite dimensional space. Denote $\mu_{\min}=\min_{x\in \Omega}\mu(x)$. Then (\ref{qp}) leads to
 $$\label{l-infty}\|u\|_{L^\infty(\Omega)}\leq \f{C}{\mu_{\min}}\|u\|_{W_0^{m,p}(\Omega)},$$
 and thus for any $1\leq q<+\infty$,
 $$\le(\int_\Omega |u|^qd\mu\ri)^{1/q}\leq \f{C}{\mu_{\min}}|\Omega|^{1/q}\|u\|_{W_0^{m,p}(\Omega)}\leq
 \f{C}{\mu_{\min}}(1+|\Omega|)\|u\|_{W_0^{m,p}(\Omega)},$$
 where $|\Omega|=\sum_{x\in\Omega}\mu(x)$ denotes the volume of $\Omega$. Therefore (\ref{Poincare}) holds. $\hfill\Box$\\

  If $V$ is a finite graph, then $W^{m,p}(V)$ can be defined as a set of all functions $u:V\ra\mathbb{R}$ under the norm
 \be\label{Wh}\|u\|_{W^{m,p}(V)}=\le(\int_V(|\nabla^m u|^p+h|u|^p)d\mu\ri)^{1/p},\ee
 where $h(x)>0$ for all $x\in V$. Then we have an obvious analog of Theorem \ref{compactlemma} as follows:
  \begin{theorem}\label{compact-V}
  Let $V$ be a finite graph. Then $W^{m,p}(V)$ is embedded in $L^q(V)$ for all $1\leq q\leq+\infty$.
  Moreover, $W^{m,p}(V)$ is pre-compact.
  \end{theorem}

 Obviously, both $W_0^{k,p}(\Omega)$  and
 $W^{k,p}(V)$ with norms  (\ref{w0-norm}) and (\ref{Wh}) respectively  are Banach spaces.
 Let $(X,\|\cdot\|)$ be a Banach space, $J: X\ra\mathbb{R}$ be a functional. We say that $J$ satisfies
 the $(PS)_c$ condition for some real number $c$, if for any sequence of functions $u_k: X\ra\mathbb{R}$
 such that $J(u_k)\ra c$ and $J^\prime(u_k)\ra 0$ as $k\ra +\infty$, there holds up to a subsequence,
 $u_k\ra u$ in $X$. To prove Theorems \ref{Theorem 1}-\ref{Theorem 8}, we need the following mountain pass theorem.
 \begin{theorem}(Ambrosetti-Rabinowitz \cite{Ambrosetti-Rabinowitz}).\label{Mountain-pass-Theorem}
  Let
 $(X,\|\cdot\|)$ be a Banach space, $J\in C^1(X,\mathbb{R})$, $e\in X$ and
 $r>0$ be such that $\|e\|>r$ and
 $$b:=\inf_{\|u\|=r}J(u)>J(0)\geq J(e).$$
 If $J$ satisfies the $(PS)_c$ condition with $c:=\inf_{\gamma\in\Gamma}\max_{t\in[0,1]}J(\gamma(t))$,
 where $$\Gamma:=\{\gamma\in C([0,1],X): \gamma(0)=0,\gamma(1)=e\},$$
 then $c$ is a critical value of $J$.
 \end{theorem}

 \section{Local existence }

 In this section, we prove Theorems \ref{Theorem 1}-\ref{Theorem 4} by applying Theorem \ref{Mountain-pass-Theorem}.
 \\

 {\it Proof of Theorem \ref{Theorem 1}}. Let $p>2$ and $\alpha<\lambda_1(\Omega)$ be fixed.
 For any $u\in W_0^{1,2}(\Omega)$, we let
 $$J(u)=\f{1}{2}\int_\Omega(|\nabla u|^2-\alpha u^2) d\mu-\f{1}{p}
 \int_\Omega (u^+)^pd\mu,$$
 where $u^+(x)=\max\{u(x),0\}$.
 It is clear that $J\in \mathcal{C}^1(W_0^{1,2}(\Omega),\mathbb{R})$. We claim that
 $J$ satisfies the $(PS)_c$ condition for any $c\in\mathbb{R}$.
 To see this, we take a sequence of functions $u_k\in W_0^{1,2}(\Omega)$ such that
 $J(u_k)\ra c$, $J^\prime(u_k)\ra 0$ as $k\ra+\infty$.
 This leads to
 \bea
 \label{p1}\f{1}{2}\int_\Omega(|\nabla u_k|^2-\alpha u_k^2)d\mu-\f{1}{p}\int_\Omega(u_k^+)^pd\mu=c+o_k(1),\\[1.5ex]
 \label{p2}\le|\int_\Omega(|\nabla u_k|^2-\alpha u_k^2)d\mu-\int_\Omega (u_k^+)^pd\mu\ri|\leq o_k(1)\|u_k\|_{W_0^{1,2}(\Omega)}.
 \eea
 Noting that $p>2$, we conclude from (\ref{p1}) and (\ref{p2}) that $u_k$ is bounded in $W_0^{1,2}(\Omega)$. Then the Sobolev
 embedding (Theorem
  \ref{compactlemma})
 implies that up to a subsequence, $u_k$ converges to some function $u$ in $W_0^{1,2}(\Omega)$. Hence the $(PS)_c$ condition holds.

  To proceed, we need to check that $J$ satisfies all conditions in the mountain pass theorem (Theorem
  \ref{Mountain-pass-Theorem}). Note that
 \be\label{j0}J(0)=0.\ee
 By Theorem \ref{compactlemma}, there exists some constant $C$ depending only on $p$ and $\Omega$ such that
 $$\le(\int_\Omega (u^+)^pd\mu\ri)^{1/p}\leq C\le(\int_\Omega|\nabla u|^2d\mu\ri)^{1/2}.$$
 Hence there holds for all $u\in W_0^{1,2}(\Omega)$
 $$J(u)\geq \f{1}{2}\|u\|_{W_0^{1,2}(\Omega)}^2-\f{C^p}{p}\|u\|_{W_0^{1,2}(\Omega)}^p.$$
 Since $p>2$, one can find some sufficiently small $r>0$ such that
 \be\label{c-2}\inf_{\|u\|_{W_0^{1,2}(\Omega)}=r}J(u)>0.\ee
  Take a function $u^*\in W_0^{1,2}(\Omega)$ satisfying $u^*>0$ in $\Omega^\circ$.
 Passing to the limit $t\ra+\infty$,  we have
 $$J(tu^*)=\f{t^2}{2}\int_\Omega(|\nabla u^*|^2-\alpha (u^*)^2)d\mu-\f{t^p}{p}\int_\Omega ({u^*}^+)^pd\mu\ra-\infty.$$
 Hence there exists some $u_0\in W_0^{1,2}(\Omega)$ such that
 \be\label{c-3}J(u_0)<0,\quad\|u_0\|_{W_0^{1,2}(\Omega)}>r.\ee
 Combining (\ref{j0}), (\ref{c-2}) and (\ref{c-3}), we conclude by Theorem \ref{Mountain-pass-Theorem} that
 $c=\inf_{\gamma\in\Gamma}\max_{t\in[0,1]}J(\gamma(t))$ is a critical value of $J$, where
 $\Gamma=\{\gamma\in C([0,1],W_0^{1,2}(\Omega)): \gamma(0)=0,\gamma(1)=u_0\}$. In particular, there exists some function
 $u\in W_0^{1,2}(\Omega)$ such that
 \be\label{eqn}-\Delta u-\alpha u=(u^+)^{p-1}\quad{\rm in}\quad \Omega^\circ.\ee
 Testing (\ref{eqn}) by $u^-=\min\{u,0\}$ and noting that $u^-u^+=u^-(u^+)^{p-1}=0$,  we have
 $$-\int_\Omega u^-\Delta ud\mu-\alpha\int_\Omega (u^-)^2d\mu=0.$$
 Since $u=u^++u^-$, the above equation leads to
 \bea\nonumber\f{\alpha}{\lambda_1(\Omega)}\int_\Omega|\nabla u^-|^2d\mu&\geq& \alpha\int_\Omega (u^-)^2d\mu\\
 &=&\nonumber
 -\int_\Omega u^-(\Delta u^++\Delta u^-)d\mu\\&=&
 \int_\Omega|\nabla u^-|^2d\mu-\int_\Omega u^-\Delta u^+d\mu.\label{g1}\eea
 Note that
 \bea\nonumber
 -\int_\Omega u^-\Delta u^+d\mu&=&-\sum_{x\in\Omega^\circ}u^-(x)\sum_{y\sim x}w_{xy}(u^+(y)-u^+(x))\\
 &=&-\sum_{x\in\Omega^\circ}\sum_{y\sim x}w_{xy}u^-(x)u^+(y)\geq 0.\label{g2}
 \eea
 Inserting (\ref{g2}) into (\ref{g1}) and recalling that $\alpha<\lambda_1(\Omega)$, we obtain
 $\int_\Omega|\nabla u^-|^2d\mu= 0$, which implies that $u^-\equiv 0$ in $\Omega$. Whence $u\geq 0$ and (\ref{eqn})
 becomes
 \be\label{eqn-1}\left\{
 \begin{array}{lll}
 -\Delta u-\alpha u=u^{p-1}\quad{\rm in}\quad \Omega^\circ\\[1.5ex]
 u\geq0\,\,\,{\rm in}\,\,\, \Omega^\circ,\,\,\, u=0\,\,\,{\rm on}\,\,\, \p\Omega.
 \end{array}
 \right.\ee
 Suppose $u(x)=0=\min_{x\in\Omega}u(x)$ for some $x\in\Omega^\circ$. If $y$ is adjacent to $x$, then
 we know from (\ref{eqn-1}) that $\Delta u(x)=0$, and thus by the definition of $\Delta$ (see (\ref{lap}) above),
 $u(y)=0$. Therefore we conclude that $u\equiv 0$ in $\Omega$, which contradicts (\ref{c-2}). Hence $u>0$ in
 $\Omega^\circ$ and this completes the proof of the theorem. $\hfill\Box$\\

 {\it Proof of Theorem \ref{Theorem 2}}. Let $p>1$ be fixed. For any $u\in W_0^{1,p}(\Omega)$, we let
 $$J_p(u)=\f{1}{p}\int_\Omega|\nabla u|^pd\mu-
 \int_\Omega F(x,u^+)d\mu,$$
 where $u^+(x)=\max\{u(x),0\}$.
 In view of $(H_4)$, there exist two constants $\lambda$ and $\delta>0$ such that $\lambda<\lambda_p(\Omega)$ and
 $$F(x,u^+)\leq \f{\lambda}{p} (u^+)^p+\f{(u^+)^{p+1}}{\delta^{p+1}}F(x,u^+).$$
 For any $u\in W_0^{1,p}(\Omega)$ with $\|u\|_{W_0^{1,p}(\Omega)}\leq 1$, we have by the Sobolev embedding
 (Theorem \ref{compactlemma}) that $\|u\|_\infty\leq C$ for some constant $C$ depending only on $p$ and $\Omega$,
 and that $$F(x,u^+)\leq \f{\lambda}{p} (u^+)^p+C{(u^+)^{p+1}}.$$
 This together with (\ref{p-eig}), the definition of $\lambda_p(\Omega)$,
 and Theorem \ref{compactlemma} leads to
 $$\int_\Omega F(x,u^+)d\mu\leq \f{\lambda}{p\lambda_p(\Omega)}\int_\Omega|\nabla u^+|^pd\mu+
 C\le(\int_\Omega|\nabla u^+|^pd\mu\ri)^{1+1/p}.$$
 Here and throughout this paper, we often denote various constants by the same $C$. Noting that
 \bna
 |\nabla u|^2(x)&=&\f{1}{2\mu(x)}\sum_{y\sim x}w_{xy}(u(y)-u(x))^2\\
 &=&\f{1}{2\mu(x)}\sum_{y\sim x}w_{xy}\le\{(u^+(y)-u^+(x))^2+(u^-(y)-u^-(x))^2\ri.\\
 &&\quad\quad\quad\le.-2u^+(y)u^-(x)-2u^+(x)u^-(y)\ri\}\\
 &\geq&|\nabla u^+|^2(x),
 \ena
  we have $|\nabla u^+|^p\leq |\nabla u|^p$. Hence if $\|u\|_{W_0^{1,p}(\Omega)}\leq 1$,
  then we obtain
 $$J_p(u)\geq \f{\lambda_p(\Omega)-\lambda}{p\lambda_p(\Omega)}\int_\Omega|\nabla u|^pd\mu
 -C\le(\int_\Omega|\nabla u|^pd\mu\ri)^{1+1/p}.$$
 Therefore
 \be\label{p-1-0}\inf_{\|u\|_{W_0^{1,p}(\Omega)}=r}J_p(u)>0\ee
 provided that $r>0$ is sufficiently small.

 By $(H_3)$, there exist two positive constants $c_1$ and $c_2$ such that
 $$F(x,u^+)\geq c_1(u^+)^{q}-c_2.$$
 Take $u_0\in W_0^{1,p}(\Omega)$ such that $u_0\geq 0$ and $u_0\not\equiv 0$. For any $t>0$, we have
 $$J_p(tu_0)\leq\f{t^p}{p}\int_\Omega|\nabla u_0|^pd\mu-c_1t^{q}\int_\Omega u_0^{q}d\mu
 -c_2|\Omega|.$$
 Since $q>p$, we conclude $J_p(tu_0)\ra -\infty$ as $t\ra+\infty$. Hence there exists some $u_1\in W_0^{1,p}(\Omega)$ satisfying
 \be\label{p-2-0}J_p(u_1)<0,\quad \|u_1\|_{W_0^{1,p}(\Omega)}>r.\ee

 We claim that $J_p$ satisfies the $(PS)_c$ condition for any $c\in\mathbb{R}$.
 To see this, we assume $J_p(u_k)\ra c$ and $J_p^\prime(u_k)\ra 0$ as $k\ra+\infty$. It follows that
 \bna
 &&\f{1}{p}\int_\Omega|\nabla u_k|^pd\mu-\int_\Omega F(x,u_k^+)d\mu=c+o_k(1)\\[1.5ex]
 &&\int_\Omega|\nabla u_k|^pd\mu-\int_\Omega u_kf(x,u_k^+)d\mu=o_k(1)\|u_k\|_{W_0^{1,p}(\Omega)}.
 \ena
 In view of $(H_3)$, we obtain from the above two equations that $u_k$ is bounded in $W^{1,p}_0(\Omega)$.
 Then the $(PS)_c$ condition follows from Theorem \ref{compactlemma}.

 Combining (\ref{p-1-0}), (\ref{p-2-0}) and the obvious fact $J_p(0)=0$, we conclude from Theorem \ref{Mountain-pass-Theorem} that
 there exists a function $u\in W_0^{1,p}(\Omega)$ such that $J_p(u)=\inf_{\gamma\in\Gamma}\max_{t\in[0,1]}J_p(\gamma(t))>0$
 and $J_p^\prime(u)=0$, where $\Gamma=\{\gamma\in C([0,1],W_0^{1,p}(\Omega)): \gamma(0)=0,\gamma(1)=u_1\}$. Hence
 there exists a nontrivial solution $u\in W_0^{1,p}(\Omega)$ to the equation
 $$-\Delta_p u=f(x,u^+)\quad{\rm in}\quad \Omega^\circ.$$
 Noting that $\Gamma(u^-,u)=\Gamma(u^-)+\Gamma(u^-,u^+)\geq |\nabla u^-|^2$, we obtain
 $$\int_\Omega|\nabla u^-|^pd\mu\leq \int_\Omega|\nabla u^-|^{p-2}\Gamma(u^-,u)d\mu=
 -\int_\Omega u^-\Delta_p ud\mu=\int_\Omega u^-f(x,u^+)dx=0.$$
 This implies that $u^-\equiv 0$ and thus $u\geq 0$. $\hfill\Box$\\

 {\it Proof of Theorem \ref{Theorem 4}}. Let $m\geq 2$ and $p>1$ be fixed. For any $u\in W_0^{m,p}(\Omega)$, we write
 $$J_{mp}(u)=\f{1}{p}\int_\Omega|\nabla^mu|^pd\mu-\int_\Omega F(x,u)d\mu.$$
  By $(A_2)$, there exists some $\lambda<\lambda_{mp}(\Omega)$ and $\delta>0$ such that for all $s\in\mathbb{R}$,
 $$F(x,s)\leq \f{\lambda}{p}|s|^p+\f{|s|^{p+1}}{\delta^{p+1}}|F(x,s)|.$$
 By Theorem \ref{compactlemma}, we have
 \bea\label{p+1}
 &&\int_\Omega|u|^{p+1}d\mu\leq C\le(\int_\Omega|\nabla^mu|^pd\mu\ri)^{1+1/p},\\
 &&\label{infty}\|u\|_{L^\infty(\Omega)}\leq C\le(\int_\Omega|\nabla^mu|^pd\mu\ri)^{1/p}.
 \eea
  For any $u\in W_0^{m,p}(\Omega)$ with $\|u\|_{W_0^{m,p}(\Omega)}\leq 1$,
 in view of (\ref{infty}), there exists some constant $C$, depending only on $m$, $p$ and $\Omega$, such that
 $\|u\|_{L^\infty(\Omega)}\leq C$ and thus $\|F(x,u)\|_{L^\infty(\Omega)}\leq C$. This together with (\ref{p+1}) gives
 \bna
 \int_\Omega F(x,u)d\mu&\leq&\f{\lambda}{p}\int_\Omega |u|^pd\mu+\f{C}{\delta^{p+1}}\int_\Omega|u|^{p+1}|F(x,u)|d\mu\\
 &\leq&\f{\lambda}{\lambda_{mp}(\Omega)p}\int_\Omega|\nabla^mu|^pd\mu+C\le(\int_\Omega|\nabla^mu|^pd\mu\ri)^{1+1/p}.
 \ena
 Hence
 $$J_{mp}(u)\geq \f{\lambda_{mp}(\Omega)-\lambda}{\lambda_{mp}(\Omega)p}\int_\Omega|\nabla^mu|^pd\mu
 -C\le(\int_\Omega|\nabla^mu|^pd\mu\ri)^{1+1/p},$$
 and thus
 \be\label{lowerb}\inf_{\|u\|_{W_0^{m,p}(\Omega)}=r}J_{mp}(u)>0\ee
 for sufficiently small $r>0$.

 By $(A_3)$, there exist two positive constants $c_1$ and $c_2$ such that
 $$F(x,s)\geq c_1|s|^{q}-c_2,\quad\forall s\in \mathbb{R}.$$
 For any fixed $u\in W_0^{m,p}(\Omega)$ with $u\not\equiv 0$, we have
 $$\label{e-0}J_{mp}(tu)\leq\f{|t|^p}{p}\int_\Omega|\nabla^mu|^pd\mu-c_1|t|^{q}\int_\Omega|u|^{q}d\mu+c_2|\Omega|\ra-\infty$$
 as $t\ra\infty$. Hence there exists some $u_2\in W_0^{m,p}(\Omega)$ such that
 \be\label{less}J_{mp}(u_2)<0,\quad \|u_2\|_{W_0^{m,p}(\Omega)}>r.\ee

 Now we claim that $J_{mp}$ satisfies the $(PS)_c$ condition for any $c\in\mathbb{R}$. For this purpose, we take
 $u_k\in W_0^{m,p}(\Omega)$ such that $J_{mp}(u_k)\ra c$ and $J^\prime_{mp}(u_k)\ra 0$ as $k\ra+\infty$. In particular,
 \bea
 \f{1}{p}\int_\Omega|\nabla^mu_k|^pd\mu-\int_\Omega F(x,u_k)d\mu=c+o_k(1),\label{p-s-1}\\[1.5ex]
 \int_\Omega|\nabla^mu_k|^pd\mu-\int_\Omega u_kf(x,u_k)d\mu=o_k(1)\|u_k\|_{W_0^{m,p}}.\label{p-s-2}
 \eea
 By $(A_3)$, we have
 \bea
 q\int_\Omega F(x,u_k)d\mu&\leq&q\int_{|u_k|\leq M}|F(x,u_k)|d\mu+q\int_{|u_k|>M}F(x,u_k)d\mu\nonumber\\
 &\leq&\int_{|u_k|> M}u_kf(x,u_k)d\mu+C\nonumber\\
 &\leq&\int_\Omega u_kf(x,u_k)d\mu+C.\label{pp1}
 \eea
 Combining (\ref{p-s-1}), (\ref{p-s-2}) and (\ref{pp1}), we conclude that $u_k$ is bounded in $W_0^{m,p}(\Omega)$.
 Then the $(PS)_c$ condition follows from Theorem \ref{compactlemma}.

 In view of (\ref{lowerb}) and (\ref{less}), applying the mountain-pass theorem as before, we finish the proof of the theorem. $\hfill\Box$

 \section{Global existence }

 In this section, using Theorem \ref{Mountain-pass-Theorem}, we prove global existence results (Theorems \ref{Theorem 5}-\ref{Theorem 8}).
 These procedures are very similar to that of Theorems \ref{Theorem 1}-\ref{Theorem 4}. We only give the outline of the proof.\\

 {\it Proof of Theorem \ref{Theorem 5}}. Let $p>2$ be fixed. For $u\in W^{1,2}(V)$, we let
 $$J_V(u)=\f{1}{2}\int_V(|\nabla u|^2+hu^2)d\mu-\f{1}{p}\int_V(u^+)^pd\mu.$$
 We first prove that $J_V$ satisfies the $(PS)_c$ condition for any $c\in\mathbb{R}$.
 To see this, we assume $J_V(u_k)\ra c$ and $J_V^\prime(u_k)\ra 0$ as $k\ra \infty$, namely
 \bea\label{psc1}&&\f{1}{2}\int_V(|\nabla u_k|^2+hu_k^2)d\mu-\f{1}{p}\int_V(u^+)^pd\mu=c+o_k(1),\\
 \label{psc2}&&\le|\int_V(|\nabla u_k|^2+hu_k^2)d\mu-\int_V(u_k^+)^pd\mu\ri|\leq o_k(1)\|u_k\|_{W^{1,2}(V)}.
 \eea
 It follows from (\ref{psc1}) and (\ref{psc2}) that $u_k$ is bounded in $W^{1,2}(V)$. By Theorem \ref{compact-V},
  $W^{1,2}(V)$ is pre-compact,
 we conclude up to a subsequence, $u_k\ra u$ in $W^{1,2}(V)$. Thus $J_V$ satisfies the $(PS)_c$ condition.

 To apply the mountain pass theorem, we check the conditions of $J_V$. Obviously
 \be\label{zer}J_V(0)=0.\ee
 Note that $h(x)>0$ for all $x\in V$. By (\ref{Wh}),
 $$\|u\|_{W^{1,2}(V)}^2=\int_V(|\nabla u|^2+hu^2)d\mu,\quad \forall u\in W^{1,2}(V).$$
 By the Sobolev embedding (Theorem \ref{compact-V}), we have for any $p>2$,
 $$\int_V (u^+)^pd\mu\leq C\|u\|_{W^{1,2}(V)}^p.$$
 Hence
 \be\label{inf}\inf_{\|u\|_{W^{1,2}(V)}=r}J_V(u)>0\ee
 for sufficiently small $r>0$.  For any fixed $\tilde{u}\geq 0$ with $\tilde{u}\not\equiv 0$, we have
 $J_V(t\tilde{u})\ra-\infty$ as
 $t\ra-\infty$. Hence there exists some $e\in W^{1,2}(V)$ such that
 \be\label{con3}J_V(e)<0,\quad \|e\|_{W^{1,2}(V)}>r.\ee
 Combining (\ref{zer}), (\ref{inf}) and (\ref{con3}) and applying the mountain-pass theorem (Theorem \ref{Mountain-pass-Theorem}),
 we find some $u\in W^{1,2}(V)\setminus\{0\}$ satisfying
 $$-\Delta u+hu=(u^+)^{p-1}\quad{\rm in}\quad V.$$
 Testing this equation by $u^-$, we have
 $$\int_V(|\nabla u^-|^2+h(u^-)^2)d\mu\leq\int_Vu^-(-\Delta u+hu)d\mu=0,$$
 which leads to $u^-\equiv 0$. Therefore $u\geq 0$ and
 $$-\Delta u+hu=u^{p-1}\quad{\rm in}\quad V.$$
 If $u(x)=0$ for some $x\in V$, then $\Delta u(x)=0$, which together with $u\geq 0$ implies that $u\equiv 0$ in $V$.
 This contradicts $u\not\equiv 0$. Therefore $u>0$ in $V$ and this completes the proof of the theorem. $\hfill\Box$\\

 {\it Proof of Theorem \ref{Theorem 6}.} Let $p>1$ be fixed. For $u\in W^{1,p}(V)$, we define
 $$J_{V,p}(u)=\f{1}{p}\int_V(|\nabla u|^p+h|u|^p)d\mu-\int_VF(x,u^+)d\mu.$$
 Write
 $$\|u\|_{W^{1,p}(V)}=\le(\int_V(|\nabla u|^p+h|u|^p)d\mu\ri)^{1/p}.$$
 By $(H_V^4)$, there exist two constants $\lambda$ and $\delta>0$ such that $\lambda<\lambda_p^h(V)$ and
 $$F(x,u^+)\leq \f{\lambda}{p} (u^+)^p+\f{(u^+)^{p+1}}{\delta^{p+1}}F(x,u^+).$$
 For any $u\in W^{1,p}(V)$ with $\|u\|_{W^{1,p}(V)}\leq 1$, we have by Theorem \ref{compact-V} that $\|u\|_\infty\leq C$,
 and that $$F(x,u^+)\leq \f{\lambda}{p} (u^+)^p+C{(u^+)^{p+1}}.$$
 This together with the definition of $\lambda_p^h(V)$ and the Sobolev embedding (Theorem \ref{compact-V})  leads to
 $$\int_V F(x,u^+)d\mu\leq \f{\lambda}{p\lambda_p^h(V)}\|u^+\|_{W^{1,p}(V)}^p
 + C\|u^+\|_{W^{1,p}(V)}^{p+1}.$$
 Then we obtain for all $u$ with $\|u\|_{W^{1,p}(V)}\leq 1$,
 $$J_{V,p}(u)\geq \f{\lambda_p^h(V)-\lambda}{p\lambda_p^h(V)}\|u\|_{W^{1,p}(V)}^p
 -C\|u\|_{W^{1,p}(V)}^{p+1}.$$
 Therefore
 \be\label{p-1}\inf_{\|u\|_{W^{1,p}(V)}=r}J_{V,p}(u)>0\ee
 provided that $r>0$ is sufficiently small.

 By the hypothesis $(H_V^3)$, there exist two positive constants $c_1$ and $c_2$ such that
 $$F(x,u^+)\geq c_1(u^+)^{q}-c_2.$$
 For some $u^*\in W^{1,p}(V)$ with $u^*\geq 0$ and $u^*\not\equiv 0$, we have for any $t>0$,
 $$J_{V,p}(tu^*)\leq\f{t^p}{p}\|u^*\|_{W^{1,p}(V)}^p-c_1t^{q}\int_V (u^*)^{q}d\mu
 -c_2|V|.$$
 Hence $J_{V,p}(tu^*)\ra -\infty$ as $t\ra+\infty$, from which one can find some $e\in W^{1,p}(V)$ satisfying
 \be\label{p-2}J_{V,p}(e)<0,\quad \|e\|_{W^{1,p}(V)}>r.\ee

 We claim that $J_{V,p}$ satisfies the $(PS)_c$ condition for any $c\in\mathbb{R}$.
 To see this, we assume $J_{V,p}(u_k)\ra c$ and $J_{V,p}^\prime(u_k)\ra 0$ as $k\ra+\infty$. It follows that
 \bna
 &&\f{1}{p}\|u_k\|_{W^{1,p}(V)}^p-\int_V F(x,u_k^+)d\mu=c+o_k(1)\\[1.5ex]
 &&\|u_k\|_{W^{1,p}(V)}^p-\int_V u_k^+f(x,u_k^+)d\mu=o_k(1)\|u_k\|_{W^{1,p}(V)}.
 \ena
 In view of $(H_V^3)$, we conclude from the above two equations that $u_k$ is bounded in $W^{1,p}(V)$.
 Then the $(PS)_c$ condition follows from Theorem \ref{compact-V}.

 Combining (\ref{p-1}), (\ref{p-2}) and $J_{V,p}(0)=0$, we conclude from Theorem \ref{Mountain-pass-Theorem} that
 there exists a nontrivial solution $u\in W^{1,p}(V)$ to the equation
 $$-\Delta_p u+h|u|^{p-2}u=f(x,u^+)\quad{\rm in}\quad V.$$
 Noting that $\Gamma(u^-,u)=\Gamma(u^-)+\Gamma(u^-,u^+)\geq |\nabla u^-|^2$, we obtain
 $$\int_V(|\nabla u^-|^p+h|u^-|^p)d\mu\leq -\int_V u^-\Delta_p ud\mu+\int_Vhu^-|u|^{p-2}ud\mu
 =\int_V u^-f(x,u^+)dx=0.$$
 This implies that $u^-\equiv 0$ and thus $u\geq 0$ in $V$. $\hfill\Box$\\

 {\it Proof of Theorem \ref{Theorem 8}.} Let $m\geq 2$ and $p>1$ be fixed. We define
 $$J_{mp}^V(u)=\f{1}{p}\int_V(|\nabla^mu|^p+h|u|^p)d\mu-\int_V F(x,u)d\mu.$$
 By $(A_V^2)$, there exists some $\lambda<\lambda_{mp}(V)$ and $\delta>0$ such that for all $x\in V$ and $s\in\mathbb{R}$,
 $$F(x,s)\leq \f{\lambda}{p}|s|^p+\f{|s|^{p+1}}{\delta^{p+1}}|F(x,s)|.$$
 For any $u\in W^{m,p}(V)$ with $\|u\|_{W^{m,p}(V)}\leq 1$, we have $\|u\|_\infty\leq C$ and thus $\|F(x,u)\|_\infty\leq C$. Hence
 \bna
 \int_V F(x,u)d\mu&\leq&\f{\lambda}{p}\int_V |u|^pd\mu+\f{C}{\delta^{p+1}}\int_V|u|^{p+1}|F(x,u)|d\mu\\
 &\leq&\f{\lambda}{\lambda_{mp}(V)p}\|u\|_{W^{m,p}(V)}^p+C\|u\|_{W^{m,p}(V)}^{p+1}.
 \ena
 It follows that
 $$J_{mp}^V(u)\geq \f{\lambda_{mp}(V)-\lambda}{\lambda_{mp}(V)p}\|u\|_{W^{m,p}(V)}^p
 -C\|u\|_{W^{m,p}(V)}^{p+1},$$
 and thus
 \be\label{lowerb10}\inf_{\|u\|_{W^{m,p}(V)}=r}J_{mp}^V(u)>0\ee
 for sufficiently small $r>0$.

 By $(A_V^3)$, there exist two positive constants $c_1$ and $c_2$ such that
 $$F(x,s)\geq c_1|s|^{q}-c_2,\quad\forall s\in \mathbb{R}.$$
 For any fixed $u\in W^{m,p}(V)$ with $u\not\equiv 0$, we have
 $$\label{e-0}J_{mp}^V(tu)\leq\f{|t|^p}{p}\int_V(|\nabla^mu|^p+h|u|^p)d\mu-
 c_1|t|^{q}\int_V|u|^{q}d\mu+c_2|V|\ra-\infty$$
 as $t\ra\infty$. Hence there exists some $e\in W^{m,p}(V)$ such that
 \be\label{less10}J_{mp}^V(e)<0,\quad \|e\|_{W^{m,p}(V)}>r.\ee

 Now we claim that $J_{mp}^V$ satisfies the $(PS)_c$ condition for any $c\in\mathbb{R}$. For this purpose, we take
 $u_k\in W^{m,p}(V)$ such that $J_{mp}^V(u_k)\ra c$ and $({J_{mp}^V})^\prime(u_k)\ra 0$ as $k\ra+\infty$. In particular,
 \bea
 \f{1}{p}\int_V(|\nabla^mu_k|^p+h|u_k|^p)d\mu-\int_V F(x,u_k)d\mu=c+o_k(1),\label{p-s-11}\\[1.5ex]
 \int_V(|\nabla^mu_k|^p+h|u_k|^p)d\mu-\int_V u_kf(x,u_k)d\mu=o_k(1)\|u_k\|_{W^{m,p}(V)}.\label{p-s-22}
 \eea
 By $(A_V^3)$, we have
 \bea
 q\int_V F(x,u_k)d\mu&\leq&q\int_{|u_k|\leq M}|F(x,u_k)|d\mu+q\int_{|u_k|>M}F(x,u_k)d\mu\nonumber\\
 &\leq&\int_{|u_k|> M}u_kf(x,u_k)d\mu+C\nonumber\\
 &\leq&\int_V u_kf(x,u_k)d\mu+C.\label{pp11}
 \eea
 Combining (\ref{p-s-11}), (\ref{p-s-22}) and (\ref{pp11}), we conclude that $u_k$ is bounded in $W^{m,p}(V)$.
 Then the $(PS)_c$ condition follows from Theorem \ref{compact-V}.

 In view of (\ref{lowerb10}), (\ref{less10}) and the fact $J_{mp}^V(0)=0$, applying the mountain-pass theorem, we finish the proof of the theorem. $\hfill\Box$\\

 {\bf Acknowledgements.} A. Grigor'yan is partly supported  by SFB 701 of the German Research Council. Y. Lin is supported by the National Science Foundation of China (Grant No.11271011). Y. Yang is supported by the National Science Foundation of China (Grant No.11171347).

\end{document}